\documentclass[leqno]{article}
\usepackage{amssymb}
\usepackage{amsmath, enumerate, vmargin}

\makeatletter
\renewcommand\section{\@startsection {section}{1}{\z@}%
 {-3.5ex \@plus -1ex \@minus -.2ex}%
 {2.3ex \@plus.2ex}%
 {\center \normalfont\large\bfseries}}
\makeatother

\setmarginsrb{3cm}{3cm}{3cm}{3cm}{1mm}{6mm}{0mm}{10mm}

\newtheorem{thm}{Theorem}[section]
\newtheorem{prop}[thm]{Proposition}
\newtheorem{cor}[thm]{Corollary}
\newtheorem{lem}[thm]{Lemma}
\newtheorem{defi}[thm]{Definition}
\newtheorem{remark}[thm]{Remark}
\newtheorem{example}[thm]{Example}
\newtheorem{pb}[thm]{Problem}

\newenvironment{rk}{\begin{remark}\rm}{\end{remark}}

\newcommand{\real}{{\mathbb R}}
\newcommand{\nat}{{\mathbb N}}

\newcommand{\com}{{\mathbb C}}

\newcommand{\E}{{\mathcal E}}
\newcommand{\F}{{\mathcal F}}

\renewcommand{\H}{{\mathcal H}}

\newcommand{\R}{{\mathcal R}}

\renewcommand{\a}{\alpha}
\renewcommand{\b}{\beta}

\newcommand{\Ga}{\Gamma}
\renewcommand{\d}{\delta}
\newcommand{\D}{\Delta}

\renewcommand{\l}{\lambda}

\newcommand{\f}{\varphi}
\renewcommand{\O}{{\Omega}}

\renewcommand{\t}{\theta}

\newcommand{\Tr}{\mbox{\rm Tr}}

\newcommand{\ot}{\otimes}
\newcommand{\op}{\oplus}
\newcommand{\8}{\infty}
\newcommand{\el}{\ell}

\newcommand{\la}{\langle}
\newcommand{\ra}{\rangle}

\newcommand{\wh}{\widehat}
\newcommand{\n}{\noindent}
\newcommand{\pf}{\noindent{\it Proof.~~}}
\newcommand{\cqd}{\hfill$\Box$}
\newcommand{\be}{\begin{eqnarray*}}
\newcommand{\ee}{\end{eqnarray*}}
\newcommand{\beq}{\begin{equation}}
\newcommand{\eeq}{\end{equation}}

\numberwithin{equation}{section}

\begin{document}


\title{Embedding of $C_q$ and $R_q$ into noncommutative
$L_p$-spaces, $1\le p<q\le 2$}
\author{Quanhua Xu}
\date{}
\maketitle


\begin{abstract}
 We prove that a quotient of subspace of
$C_p\op_pR_p$ ($1\le p<2$) embeds completely isomorphically into a
noncommutative $L_p$-space, where $C_p$ and $R_p$ are respectively
the $p$-column and $p$-row Hilbertian operator spaces. We also
represent $C_q$ and $R_q$ ($p<q\le2$) as  quotients of subspaces
of $C_p\op_pR_p$. Consequently, $C_q$ and $R_q$ embed completely
isomorphically into a noncommutative $L_p(M)$. We further show
that the underlying von Neumann algebra $M$ cannot be semifinite.
\end{abstract}

 \setcounter{section}{-1}

\makeatletter
\renewcommand{\@makefntext}[1]{#1}
\makeatother \footnotetext{\noindent Laboratoire de
Math\'{e}matiques,
 Universit\'{e} de Franche-Comt\'{e},
 25030 Besan\c con, cedex - France\\
 qx@math.univ-fcomte.fr
 }

\section{Introduction}

In the recent remarkable work \cite{ju-OH}, Junge proved that
Pisier's space $OH$ (the operator Hilbert space) embeds completely
isomorphically into the predual of a QWEP von Neumann algebra. He
further showed in \cite{ju-OH hyperfinite} that the von Neumann
algebra in question can be chosen to be hyperfinite. Junge's proof
consists of two steps. The first one is to represent $OH$ (in a
precise way) as a quotient of subspace of $C\op R$, where $C$ and
$R$ are the usual column and row Hilbertian operator spaces. (By a
quotient of subspace of $X$ we mean a quotient of a subspace $Y$
of $X$.) His tool for this representation theorem is a formula due
to Pusz and Woronowicz on the square root of the product of two
commuting positive operators on a Hilbert space (cf. [PW]). The
second step is a kind of Khintchine type inequality involving a
three term K-functional. Junge obtained such an inequality via
reduced free product by using Speicher's central limit theorem
(cf. \cite{speicher}).

Junge's arguments were simplified by Pisier \cite{pis-Rp}. The
fact that $OH$ is a quotient of subspace of $C\op R$ has been
already known before \cite{ju-OH} (see \cite [Exercise
7.9]{pis-intro}). Pisier's approach to this is via the complex
interpolation. Junge's second step above was then replaced in
\cite{pis-Rp} by a Khintchine type inequality for the generalized
circular systems obtained in \cite{pisshlyak}.

The fact that the complex interpolation can be used instead of
Pusz-Woronowicz's formula is not surprising since this formula is
intimately related to interpolation. Indeed, it can be shown that
Pusz-Woronowicz's formula is equivalent to the real interpolation
of weighted $L_2$-spaces. This will be discussed in details in the
forthcoming work \cite{jx-quantum} (see also section~3 below for a
brief discussion).

\medskip

The purpose of the present paper is to consider the same embedding
problem for the  $p$-column and  $p$-row spaces $C_p$ and $R_p$.
Recall that $C_p$ (resp. $R_p$) is the subspace of the Schatten
class $S_p$ consisting of matrices whose all entries but those in
the first column (resp. row) vanish. Both $C_p$ and $R_p$ are
isometric to $\ell_2$ as Banach spaces. Note that $C_2$ and $R_2$
coincide with $OH$ (completely isometrically), and $C_\8$ and
$R_\8$ are just $C$ and $R$ respectively.

We will prove that if $1\le p<q\le2$, then both $C_q$ and $R_q$
embed completely isomorphically into a noncommutative $L_p$-space.
An intermediate result is that $C_q$ and $R_q$ are quotients of
subspace of $C_p\op_pR_p$. We will present in section~3 two proofs
for this representation result. The first one is via the complex
interpolation as in \cite{pis-Rp}, whereas the second is based on
the real interpolation.

Besides the representation of $C_q$ as a quotient of subspace of
$C_p\op_pR_p$, another ingredient of the embedding theorem above
is the Khintchine type inequality for generalized circular systems
in noncommutative $L_p$-spaces from \cite{xu-gro}. We will recall
in section~4 this inequality, and also give its Fermionic analogue
(at least for $p>1$).

More generally, we will show that any quotient of subspace of
$C_p\op_pR_p$ ($1\le p<2$) embeds completely isomorphically into a
noncommutative $L_p$. This will be done by the Khintchine
inequality quoted previously and a decomposition result for
subspaces of $C_p\op_pR_p$ (see Proposition~\ref{graph}). It was
proved in \cite{xu-gro} that if an operator space $X$ and its dual
embed completely isomorphically into a noncommutative $L_p$
($1<p<2$), then $X$ is completely isomorphic to a quotient of
subspace of $H_p^c\op_p K_p^r$ for two Hilbert spaces $H$ and $K$,
where $H_p^c$ (resp. $K_p^r$) means that $H$ (resp. $K$) is
equipped with the $p$-column (resp. $p$-row) operator space
structure. Therefore, we obtain that for $1<p<2$ an operator space
$X$ is completely isomorphic to a quotient of subspace of
$H_p^c\op_p K_p^r$ iff both $X$ and $X^*$ embed completely
isomorphically into a noncommutative $L_p$.

To end this introduction let us mention that all algebras
constructed in \cite{ju-OH, ju-OH hyperfinite, pis-Rp} are of type
III. This is not a hazard for Pisier \cite{pis-III} proved that
$OH$ cannot embed completely isomorphically into the predual of a
semifinite von Neumann algebra. We will also extend this result by
Pisier to our general setting: $C_q$ cannot embed completely
isomorphically into  a semifinite $L_p$ ($1\le p<q\le2$).

\section{Preliminaries}

In this section we recall briefly the natural operator space
structure on noncommutative $L_p$-spaces and some elementary
properties of column and row spaces. We use standard notions and
notation from operator space theory (cf. \cite{er-book,
pis-intro}).

\medskip

It is well known that there are several equivalent constructions
of noncommutative $L_p$-spaces associated with a von Neumann
algebra. All these constructions give the same spaces (up to
isometry). In this paper we will use Haagerup noncommutative
$L_p$-spaces constructed in \cite{haag-Lp}. \cite{terp-Lp} gives
an excellent exposition of the subject. Let $M$ be a von Neumann
algebra and $1\le p\le\8$. As usual, $L_p(M)$ denotes the Haagerup
$L_p$-space associated with $M$. Recall that $L_\8(M)$ is just $M$
itself, and $L_1(M)$ can be identified with the predual $M_*$ of
$M$. More generally, if $1\le p<\8$ and $1/p+1/p'=1$, the dual of
$L_p(M)$ is $L_{p'}(M)$. When $M$ is semifinite, we will always
consider $L_p(M)$ as the usual $L_p$-space constructed from a
normal semifinite faithful trace.

Now we  describe the natural operator space structure on $L_p(M)$
as introduced in \cite{pis-ast, pis-intro}. For $p=\infty$,
$L_p(M)=M$ has its natural operator space structure as a von
Neumann algebra. For $p=1$, the natural operator space structure
on $L_1(M)=M_*$ is induced by the natural embedding of $M_*$ into
its bidual $M^*$, the latter being the standard dual of $M$. In
fact, as explained in \cite[section~7]{pis-intro}, it is more
convenient to consider $L_1(M)$ as the predual of the opposite von
Neumann algebra $M^{\rm op}$, which is isometric (but in general
not completely isomorphic) to $M$, and to equip $L_1(M)$ with the
operator space structure inherited from $(M^{\rm op})^*$. The main
reason for this choice is that the duality between $S_\8^n$ and
$S_1^n$ in operator space theory is the parallel duality (instead
of the usual tracial duality). Note that this choice further
insures that the equality $L_1(M_n\ot M)=S_1^n\wh \ot L_1(M)$
(operator space projective tensor product) holds true. Finally,
the operator space structure on $L_p(M)$ is obtained by complex
interpolation, using the well known interpretation of $L_p(M)$ as
the interpolation space $(M,\; L_1(M))_{1/p}$ (see
\cite{terp-int}). If $M$ admits a normal faithful state $\f$, we
can also use Kosaki's interpolation \cite{kos-int}.

\medskip

In particular, the Schatten classes $S_p$ are equipped with their
natural operator space structure. We will need vector-valued
Schatten classes as introduced in \cite{pis-ast}. Let $E$ be an
operator space. $S_1[E]$ is defined as the operator space
projective tensor product $S_1 \widehat \otimes E$. Then for any
$1<p<\8$, $S_p[E]$ is defined by interpolation:
 $$S_p[E]=\big(S_\8[E],\ S_1[E]\big)_{1/p}\;.$$

When $E$ is a subspace of a noncommutative $L_p(M)$, the norm of
$S_p[E]$ admits a very simple description.  Indeed, $S_p[L_p(M)]$
coincides (completely isometrically) with
$L_p(B(\ell_2)\overline{\ot}M)$ ($p<\8$).  If $E\subset L_p(M)$,
then $S_p[E]$ is the closure in $L_p(B(\ell_2)\overline{\otimes}
M)$ of the algebraic tensor product $S_p\ot E$.

By reiteration, for any $1\le p_0, p_1\le\8$ and $0<\t<1$
 $$\big(S_{p_0}[E],\ S_{p_1}[E]\big)_{\t}=S_p[E],$$
where $\frac{1}p=\frac{1-\t}{p_0}+\frac{\t}{p_1}$. More generally,
given a compatible couple $(E_0, E_1)$ of operator spaces we have
 \beq\label{vector valued interplation}
 \big(S_{p_0}[E_0],\ S_{p_1}[E_1]\big)_{\t}
 =S_p[(E_0,\ E_1)_\t].
 \eeq

The main concern of this paper is on the $L_p$-space counterparts
of the usual row and column spaces $R$ and $C$. Let $C_p$ (resp.
$R_p$) denote the subspace of $S_p$ consisting of matrices whose
all entries but those in the first column (resp. row) vanish.  It
is clear that $C_p$ and $R_p$ are completely 1-complemented
subspaces of $S_p$. We have the following completely isometric
identifications:
 \beq\label{duality Cp-Rp}
 (C_p)^*\cong C_{p'}\cong R_p\quad\mbox{and}\quad
 (R_p)^*\cong R_{p'}\cong C_p\,,\quad\forall\;1\le p\le\8.
 \eeq

 $C_p$ and $R_p$ can  be also defined via interpolation from $C$
and $R$. We view $(C, R)$ as a compatible couple by identifying
both of them with $\ell_2$ (at the Banach space level), i.e. by
identifying the canonical bases $(e_{k,1})$ of $C_p$ and
$(e_{1,k})$ of $R_p$ with $(e_k)$ of $\ell_2$. Then
 $$C_p=(C,\ R)_{1/p}=(C_\8,\ C_1)_{1/p}\quad\mbox{and}\quad
 R_p=(R,\ C)_{1/p}=(R_\8,\ R_1)_{1/p}.$$
By reiteration, for any $1\le p_0, p_1\le\8$ and $0<\t<1$
 \beq\label{interplation between Cp}
 C_p=(C_{p_0},\ C_{p_1})_{\t}\quad\mbox{and}\quad
 R_p=(R_{p_0},\ R_{p_1})_{\t},
 \eeq
where $\frac{1}{p}=\frac{1-\t}{p_0}+\frac{\t}{p_1}$. Like $C$ and
$R$, $C_p$ and $R_p$ are 1-homogenous 1-Hilbertian operator
spaces. We refer to \cite{pis-ast, pis-OH} for the proofs of all
these elementary facts.

More generally, given a Hilbert space $H$ and $1\le p\le\8$ we
denote by $H_p^c$ (resp. $H_p^r$) the Schatten $p$-class
$S_p(\com, H)$ (resp. $S_p(\bar H, \com)$) equipped with its
natural operator space structure. $H_p^c$ (resp. $H_p^r$) can be
naturally viewed as the column (resp. row) subspace of the
Schatten class $S_p(H)$ (resp. $S_p(\bar H)$). When $H$ is
separable and infinite dimensional, $H_p^c$ and $H_p^r$ are
respectively $C_p$ and $R_p$ above.  All previous properties for
$C_p$ and $R_p$ hold for $H_p^c$ and $H_p^r$ too. $H_p^c$ and
$H_p^r$ are respectively the $p$-column and $p$-row spaces
associated with $H$.

Now let $E$ be an operator space. We denote by $C_p[E]$ (resp.
$R_p[E]$) the closure of $C_p\ot E$ (resp. $R_p\ot E$) in
$S_p[E]$. Again, $C_p[E]$ and  $R_p[E]$ are completely
1-complemented subspaces of $S_p[E]$. If $E$ is a subspace of a
noncommutative $L_p(M)$, the norm of $C_p[E]$ is easy to be
determined as follows. For any finite sequence $(x_k)\subset E$
 $$\|\sum_kx_k\ot e_k\|_{C_p[E]}
 =\|(\sum_kx_k^*x_k)^{1/2}\|_{L_p(M)}\;,$$
where $(e_k)$ denotes the canonical basis of $C_p$. More
generally, if $a_k\in C_p$, then
 \beq\label{description of Cp(Lp)}
 \|\sum_kx_k\ot a_k\|_{C_p[E]}
 =\|(\sum_{k,j}\la a_k, a_j\ra x_k^*x_j)^{1/2}\|_{L_p(M)}\;,
 \eeq
where $\la\; ,\;\ra$ denotes the scalar product in $C_p$. (In
terms of matrix product, $\la a_k, a_j\ra=a_k^*a_j$.) We have also
a similar description for the norm of $R_p[E]$. The same formulas
hold for more general $H_p^c[E]$ and $H_p^r[E]$.

\section{A decomposition}

Let $T: X\to Y$ be a linear map between two Banach spaces. Recall
that the graph (or domain) of $T$ is the subspace of $X\op Y$
given by
 $$G(T)=\{(x, Tx)\ :\ x\in {\rm Dom}(T)\},$$
where ${\rm Dom}(T)$ stands for the definition domain of $T$.

\begin{prop}\label{graph}
 Let $X_0, X_1$ be 1-homogeneous 1-Hilbertian operator
spaces. Let $X= X_0\oplus_2 X_1$ and $Y\subset X$ be a closed
subspace. Then there are closed subspaces $ Y_0, Z_0\subset X_0,\
Y_1, Z_1\subset X_1$ and an injective closed densely defined
operator $T$ from $Z_0$ to $Z_1$ with dense range such that
 $$Y=Y_0  \op_2 Y_1  \op_2 G(T)\quad
 \mbox{with completely equivalent norms},$$
where $G(T)$ is the graph of $T$ regarded as a subspace of $Z_0
\oplus_2 Z_1\subset X_0\oplus_2 X_1 $. Moreover, the relevant
constants in the above equivalence are universal.
 \end{prop}

\pf Let $\la\;,\;\ra_i$ denote the scalar product  of $X_i$, and
let $P_i$  be the natural projection from $X$ onto $X_i$ ($i=0,
1$). Since the map $(x,y)\mapsto \la P_0(x),\; P_0(y)\ra_0$ is a
contractive sesquilinear form on $Y$, there is  a contraction $A:
Y\to Y$ such that
 \beq\label{graph1}
 \la x,\; A(y)\ra =\la P_0(x),\;
 P_0(y)\ra_0,\quad  x, y\in Y.
 \eeq
It is clear that $A\ge 0$ and
 \beq\label{graph2}
 \la x,\; (1-A)(y)\ra =\la P_1(x),\;
 P_1(y)\ra_1,\quad  x, y\in Y.
 \eeq
Therefore,
 $$
 x\in{\rm ker}\,A \Leftrightarrow P_0(x)=0
 \quad \mbox{and}\quad
 x\in{\rm ker}\,(1-A) \Leftrightarrow P_1(x)=0.
 $$
Set
 $$Y_0={\rm ker}\,(1-A)\quad \mbox{and}\quad
  Y_1={\rm ker}\,A.$$
Then it follows that
 $$P_i\big|_{Y_i}=A\big|_{Y_i},\quad i=0, 1.$$
Consequently,
 $$Y_i\subset X_i,\quad i=0, 1.$$
Let $Z$ be the orthogonal complement of $Y_0  \oplus_2 Y_1$ in
$Y$. Then we have the following orthogonal decomposition of $Y$ at
the Banach space level:
 \beq\label{graph3}
 Y=Y_0 \oplus_2 Y_1\oplus_2 Z.
 \eeq
We will show that this decomposition also holds in the category of
operator spaces with completely equivalent norms.

Let $y_0\in Y_0$ and $z\in Z$. Then by (\ref{graph1}) and the
definition of $Y_0$
 $$\la y_0,\; P_0(z)\ra_0=\la P_0(y_0),\; P_0(z)\ra_0
 =\la y_0,\; A(z)\ra=\la A(y_0),\; z\ra=\la y_0,\; z\ra=0.$$
Therefore, $Y_0$ and $P_0(Z)$ are orthogonal subspaces of $X_0$.
Let $Q_0$ be the orthogonal projection from $X_0$ onto $Y_0$.
Since $X_0$ is homogeneous, $Q_0$ is completely contractive. Put
$R_0=Q_0P_0\big|_{Y}$. Then $R_0$ is a completely contractive
projection from $Y$ onto $Y_0$. In fact, it is easy to see that
$R_0$ is exactly the orthogonal projection onto $Y_0$. Indeed, we
have $R_0=0$ on $Y_1$, and as well as  on $Z$ for
$Q_0\big|_{P_0(Z)}=0$. Thus $R_0=0$ on the orthogonal complement
of $Y_0$. Similarly, the orthogonal projection $R_1$ from $Y$ on
$Y_1$ is also completely contractive. It follows that the
decomposition in (\ref{graph3}) holds at the operator space level
with completely equivalent norms.

It remains to express the third space $Z$ as a graph. To that end
let us observe that both $A$ and $1-A$ are injective operators on
$Z$ with dense range. Consequently, the inverses of $A$ and $1-A$
are densely defined operators on $Z$ too.

On the other hand, since $A$ and $P_0\big|_Y$ (resp. $1-A$ and
$P_1\big|_Y$) have the same kernel, we deduce that any $z\in Z$ is
uniquely determined by $P_0(z)$ as well as by $P_1(z)$. Therefore,
$Z$ must be the graph of a closed densely defined operator, i.e.
the graph of the map $P_0(z)\mapsto P_1(z)$. Let us explain this
elementary fact precisely in the present situation. Let $Z_i$ be
the closure of $P_i(Z)$ in $X_i$. By (\ref{graph1}) and
(\ref{graph2}), there are unitaries $u_i$ from $Z_i$ to $Z$ such
that
 $$A^{1/2}\big|_{Z}= u_0P_0\big|_{Z} \quad \mbox{and}\quad
  (1-A)^{1/2}\big|_{Z}= u_1P_1\big|_{Z}\;.$$
Therefore, for any $z\in Z$ we have
 $$P_1(z)=u_1^*(1-A)^{1/2}(z)
 =u_1^*(1-A)^{1/2}A^{-1/2}u_0\big(P_0(z)\big).$$
Thus
 $$T=u_1^*(1-A)^{1/2}A^{-1/2}u_0$$
is a closed densely defined operator from $Z_0$ to $Z_1$ whose
domain coincides with $Z$. Note that $T$ is invertible and its
inverse is also densely defined. This completes the proof.\cqd

\medskip

\begin{rk}\label{graphbis}
 Let $T=w\Delta$ be the polar decomposition of $T$. Then $w$
is a unitary and $\Delta$ is a positive self-adjoint non-singular
closed densely defined operator on $Z_0$.  By the homogeneity of
$X_1$ and replacing $Z_1$ by $Z_0$ (but still equipped with the
operator space structure of $X_1$), we can assume $T=\Delta$ is
positive. Let us note, however, that with the previous precise
form of $T$ we can simply take $(1-A)^{1/2}A^{-1/2}$ as $\Delta$.
\end{rk}

Applying the previous proposition to row and column spaces, we get
the following

\begin{cor}\label{graphter}
 Let $H, K$ be Hilbert spaces and $1\le p\le\8$. Let $E$  be
a closed subspace of $H^c_p\op_p K^r_p$. Then there are closed
subspaces $ H_0, H_1\subset H,\ K_0, K_1\subset K$ and an
injective closed densely defined operator $D$ from $H_1$ to $K_1$
with dense range such that
 $$E=(H_0)^c_p  \op_p (K_0)^r_p  \op_p G(D)\quad
 \mbox{with completely equivalent norms},$$
where $G(D)$ is the graph of $D$ regarded as a subspace of
$(H_1)^c_p\op_p (K_1)^r_p$. Moreover, the relevant constants in
the above equivalence are universal.
 \end{cor}

\n{\bf Remark.} In the situation of the above corollary, it
suffices to consider the case of $p\le 2$ or $p\ge2$ for we have
$$H_{p}^c\cong H_{p'}^r\quad
 \mbox{completely isometrically}\,.$$

\section{$C_q$ as a quotient of subspace of $C_p\op_p R_p$}

It is known that $OH$ can be represented as a quotient of subspace
of $C\op R$ (see \cite[Exercise 7.9]{pis-intro} and
\cite{pis-Rp}). Pisier's proof is via complex interpolation by
using his interpolation formula $OH=(C,\ R)_{1/2}$ obtained in
\cite{pis-OH}. The aim of this section is to show a similar result
for $C_q$ relative to the couple $(C_p,\ R_p)$ for any $q$ between
$p$ and $p'$ (recalling that $p'$ denotes the index conjugate to
$p$).

\begin{thm}\label{Cq as a qs of Cp+Rp}
 Let $1\le p\le\8$ and $\min(p, p')< q< \max(p, p')$. Then
$C_q$ and $R_q$ are completely isometric to  quotients of subspace
of $C_p\op_p R_p\,$.
 \end{thm}

We will present two proofs for the theorem above. The first is
based on the complex interpolation like in \cite{pis-Rp}. The
second one is similar to Junge's approach to the representation of
$OH$ as a quotient of subspace of $C\op_\8 R$ (cf. \cite{ju-OH});
but instead of Pusz-Woronowicz's formula, we will use directly the
real interpolation. Both approaches have their advantages. The
first one is more transparent and slightly simpler. Moreover, it
shows that $C_q$ is completely \underbar{isometric} to a quotient
of subspace of $C_p\op_p R_p$. The major advantage of the second
proof is that it gives the precise form of the quotient of
subspace (or the two densities involved in this quotient) of
$C_p\op_p R_p$ which is completely isomorphic to $C_q$. This
precise form is crucial for many applications.

\medskip

We will need the following lemma from \cite{xu-gro}, which is a
generalization of \cite[Theorem~8.4]{pis-OH}.  See also
\cite{xu-description} for a more general result.

\begin{lem}\label{description of norm in SpCq}
 Let $2\le p\le\8,\; 0<\t<1,\; \frac{1}{q}
 =\frac{1-\t}{p}+\frac{\t}{p'}$
and $r$ be the conjugate index of $p/2$. Let $(e_k)$ denote the
canonical basis of $C_q$. Then for any finite sequence
$(x_k)\subset S_p$
 $$\|\sum_kx_k\ot e_k\|_{S_p[C_q]}
 =\sup\Big\{\big(\sum_k\|\a x_k\b\|_2^2\big)^{1/2}\Big\},$$
where the supremum runs over all $\a$ and $\b$ in the unit ball of
$S_{2r\t^{-1}}$ and $S_{2r(1-\t)^{-1}}$, respectively. Moreover,
the supremum can be restricted to all $\a$ and $\b$ in the
positive parts of these unit balls.
 \end{lem}

\n{\em Proof of Theorem~\ref{Cq as a qs of Cp+Rp}.} This proof is
similar to the corresponding proof on the representation of $OH$
as a quotient of subspace of $C\op R$ in \cite{pis-Rp}. Since
$R_q\cong C_{q'}$, it suffices to prove the assertion for $C_q$.
By symmetry, we may assume $p\ge 2$. We will use the following
interpolation formula (see (\ref{duality Cp-Rp}) and
(\ref{interplation between Cp}))
 \beq\label{Cq as a qs of Cp+Rp1}
 C_q=(C_p,\ R_p)_{\t},
 \eeq
where $\t$ is determined by
$\frac{1}{q}=\frac{1-\t}{p}+\frac{\t}{p'}$.   To proceed further
to the proof, we need to recall some elementary facts on complex
interpolation. Let $(X_0,\ X_1)$ be a compatible couple of Banach
spaces. Set $S=\{z\in\com\ :\ 0\le{\rm Re}(z)\le 1\}$. Let $\mu$
be the harmonic measure in $S$ at the point $\t$. Note that $\mu$
can be written as $\mu=(1-\t)\mu_0 + \t\mu_1$, where $\mu_0$ and
$\mu_1$ are probability measures respectively on the left boundary
$\partial_0$ and the  right boundary $\partial_1$ of $S$. Let
${\cal F}(X_0, X_1)$ be the family of all functions $f$ from $S$
to $X_0+X_1$ satisfying the following conditions:
 \begin{enumerate}[i)]
 \item $f$ is analytic in the interior of $S$;
 \item $f\big|_{\partial_j}\in L_2(X_j; \mu_j),\; j=0,1$.
 \end{enumerate}
 ${\cal F}(X_0, X_1)$ is equipped with the norm
 \beq\label{Cq as a qs of Cp+Rp2}
 \|f\|=\big[(1-\t)\|f\big|_{\partial_0}\|_{L_2(X_0; \mu_0)}^p
 +\t \|f\big|_{\partial_1}\|_{L_2(X_1; \mu_1)}^p\big]^{\frac{1}p}
 \eeq
(with the obvious change for $p=\8$). Then the complex
interpolation space $(X_0,\ X_1)_{\t}$ can be defined as the space
of all $x\in X_0+X_1$ such that there is $f\in {\F}(X_0, X_1)$
such that $f(\t)=x$, equipped with the quotient norm
 $$\|x\|_\t=\inf\big\{\|f\|_{{\F}(X_0, X_1)}\ :\ f(\t)=x,\;
 f\in {\F}(X_0, X_1)\big\}.$$
We refer to \cite{bl} for more details. We should call the
reader's attention to the fact that the $\ell_p$-norm in (\ref{Cq
as a qs of Cp+Rp2}) is usually replaced by the $\ell_2$-norm (to
be  consistent with the $L_2$-norm). However, at the end we get
the same interpolation norm $\|\  \|_\t$. (In fact, the indices
$2,\; p$ in (\ref{Cq as a qs of Cp+Rp2}) can be replaced by any
two indices in $[1, \8]$ without changing the norm $\|\ \|_\t$.)
Note that ${\F}(X_0, X_1)$ is a subspace of
 $${\E}(X_0, X_1)=(1-\t)^{\frac{1}p}L_2(X_0; \mu_0)\op_p
 \t^{\frac{1}p}L_2(X_1;\mu_1).$$
Using the Poisson integration, we easily see that ${\E}(X_0, X_1)$
is exactly the space of all  functions from $S$ into $X_0+X_1$
harmonic inside $S$ and verifying the condition ii) above. Then
${\F}(X_0, X_1)$ is the subspace of ${\E}(X_0, X_1)$ consisting of
analytic functions.

\medskip

Now return back to the formula (\ref{Cq as a qs of Cp+Rp1}). In
view of the above discussion set
 $$H=(1-\t)^{\frac{1}p}L_2(\ell_2; \mu_0)\quad\mbox{and}
 \quad K=\t^{\frac{1}p}L_2(\ell_2; \mu_1).$$
$H$ and $K$ are infinite dimensional separable Hilbert spaces. Let
 $$E=H^c_p\op_p K^r_p.$$
As mentioned above, $E$ can be considered as the space of
$\ell_2$-valued harmonic functions in $S$  with square integrable
boundary values. Let $F\subset E$ be the subspace of analytic
functions. Thus $F$ is a subspace of $C_p\op_p R_p$. Consider the
map $Q: F\to C_q$ defined by $Q(f)=f(\t)$. We are going to show
that $Q$ is a completely isometric surjection. By virtue of
\cite[Lemma~1.7]{pis-ast}, this is equivalent to show that
$I_{S_p}\ot Q$ extends to a quotient map from $S_p[F]$ to
$S_p[C_q]$.

First we prove that $I_{S_p}\ot Q$ extends to a contraction from
$S_p[F]$ to $S_p[C_q]$ (which amounts to saying $\|Q\|_{cb}\le 1$
). Let $f\in S_p[F]$ be of norm $\le1$ (assuming that $f$ is a
finite matrix). Let $x=f(\t)$. We have to show
$\|x\|_{S_p[C_q]}\le 1$. To this end we use Lemma \ref{description
of norm in SpCq}. Let $f_n$ be the n$^{\rm th}$ coordinate of $f$.
Set $x_n=f_n(\t)$. Let $\a$ (resp. $\b$) be a positive unit
element in $S_{2r\t^{-1}}$ (resp. $S_{2r(1-\t)^{-1}}$). Define for
$z\in S$
 $$g(z)=\sum_ng_n(z)\ot e_n\quad \mbox{with}\quad
 g_n(z)=\a^{\frac{z}{\t}}\,f_n(z)\,\b^{\frac{1-z}{1-\t}}.$$
Then $g_n$ is an analytic function in $S$ with values in $S_2$ and
 $$x_n=g_n(\t)=(1-\t)\int_{\partial_0}\,g_n(z)\,
 d\mu_0(z)+\t\int_{\partial_1}\,g_n(z)\,d\mu_1(z).$$
By the H\"older inequality, we have
 \be
 \int_{\partial_0}\sum_n\|g_n(z)\|_2^2\,d\mu_0(z)
 &\le&\int_{\partial_0}\sum_n
 \big\|f_n(z)\b^{\frac{1}{1-\t}}\big\|_2^2\,d\mu_0(z)\\
 &=&\int_{\partial_0}\sum_n{\rm Tr}\big(\b^{\frac{1}{1-\t}}
 f_n(z)^*f_n(z)\b^{\frac{1}{1-\t}}\big)\,d\mu_0(z)\\
 &\le&\Big\|\sum_n\int_{\partial_0}
  f_n(z)^*f_n(z)d\mu_0(z)\Big\|_{p/2}.
 \ee
Similarly,
 $$\int_{\partial_1}\sum_n\|g_n(z)\|_2^2\,d\mu_1(z)
 \le\Big\|\sum_n\int_{\partial_1}f_n(z)
 f_n(z)^*d\mu_1(z)\Big\|_{p/2}.$$
On the other hand, noting
 $$S_p[E]=S_p[H^c_p]\op_p S_p[K^r_p]=
 H^c_p[S_p]\op_p H^r_p[S_p]$$
and using (\ref{description of Cp(Lp)}), we have
 \be
 \|f\|_{S_p[F]}^p
 &=&(1-\t)\Big\|\sum_n\int_{\partial_0}
 f_n(z)^*f_n(z)d\mu_0(z)\Big\|_{p/2}^{p/2}\\
 &&\hskip 0.5cm+\t\,\Big\|\sum_n\int_{\partial_1}
 f_n(z)f_n(z)^*d\mu_1(z)\Big\|_{p/2}^{p/2}.
 \ee
Combining the preceding inequalities with the classical Jensen
inequality, we get (recalling that $p\ge 2$)
 \be
 \big(\sum_n\|\a x_n\b\|_2^2\big)^{1/2}
 &=&\big(\sum_n\|g_n(\t)\|_2^2\big)^{1/2}\\
 &\le&
 \big[(1-\t)\int_{\partial_0}\sum_n\|g_n(z)\|_2^2\,d\mu_0(z)
 +
 \t\int_{\partial_1}\sum_n\|g_n(z)\|_2^2\,d\mu_1(z)\big]^{1/2}\\
 &\le& \|f\|_{S_p[F]}\le 1.
 \ee
Therefore, by Lemma~\ref{description of norm in SpCq}, we deduce
$\|x\|_{S_p[C_q]}\le 1$, and so $Q$ is completely contractive.

\medskip

Now let $x\in S_p[C_q]$ with $\|x\|_{S_p[C_q]}< 1$. We have to
find $f\in S_p[F]$ such that $f(\t)=x$ and $\|f\|_{S_p[F]}<1$.
This is easy by interpolation. Indeed, by (\ref{Cq as a qs of
Cp+Rp1}) and (\ref{vector valued interplation})
 $$S_p[C_q]=\big(S_p[C_p],\ S_p[R_p]\big)_{\t}.$$
Thus there is $f\in {\F}(S_p[C_p], S_p[R_p])$ such that $f(\t)=x$
and
 $\|f\|_{{\F}(S_p[C_p], S_p[R_p])}<1$.
Using the notation in the first part of the proof, we have
 $$\Big\|\sum_n\int_{\partial_0}
  f_n(z)^*f_n(z)d\mu_0(z)\Big\|_{p/2}
 \le \big\|f\big|_{\partial_0}\big\|_{L_2(\mu_0; S_p[C_p])}^2$$
and a similar inequality for the right endpoint. It then follows
that $f\in S_p[F]$ and
 $$\|f\|_{S_p[F]}\le\|f\|_{{\F}(S_p[C_p], S_p[R_p])}<1.$$
Therefore, $Q$ is a complete quotient map, and so $C_q$ is
completely isometric to the quotient space $F/\ker(Q)$.\cqd

\medskip

We turn to the second approach to Theorem \ref{Cq as a qs of
Cp+Rp} via the real interpolation. Below we will be rather brief
and refer to \cite{jx-quantum} for more details and for more
``concrete'' representations of $C_q$ as a quotient of subspace of
$C_p\op_pR_p$.

We start by recalling the real interpolation method (see \cite{bl}
for more details). Let $(X_0, X_1)$ be a compatible couple of
Banach spaces. Let $X_0 +_p X_1$ denote the quotient of $X_0 \op_p
X_1$ by the subspace $\{(x_0, x_1): x_0+x_1=0\}$. Similarly, we
denote by $X_0 \cap_p X_1$ the diagonal subspace of $X_0 \op_p
X_1$. Then $(X_0, X_1)_{\t, p;K}$ is defined as the subspace of
 $$L_p(X_0; t^{-\t})+_p\, L_p(X_1; t^{1-\t})$$
consisting of all constant functions. Here for a given space $X$
and $\a\in\real$ we have denoted by $L_p(X; t^\a)$ the $X$-valued
$L_p$-space over $\real_+$ with respect to the weighted measure
$t^{p\a}dt/t$.  $(X_0, X_1)_{\t, p;K}$ can be also described by
the J-method.  Let $(X_0, X_1)_{\t, p;J}$ be the quotient space of
 $$L_p(X_0; t^{-\t})\cap_p\, L_p(X_1; t^{1-\t})$$
by the subspace of all $f$ such that $\int_0^\8 f(t)dt/t=0$. Then
$(X_0, X_1)_{\t, p;K}=(X_0, X_1)_{\t, p;J}$ isomorphically with
relevant constants depending only on $\t$ and $p$.

\medskip

Now let
 $$K_{\t,p}=(L_2(\el_2; t^{-\t}))_p^c
 \,+_p\, (L_2(\el_2; t^{1-\t}))_p^r$$
(recall that $(L_2(\el_2; t^{\a}))_p^c$ means the Hilbert space
$L_2(\el_2; t^{\a})$ equipped with the $p$-column space
structure). We equip $K_{\t,p}$ with the quotient operator space
structure given by that of  $(L_2(\el_2; t^{-\t}))_p^c
 \,\op_p\, (L_2(\el_2; t^{1-\t}))_p^r$.
Let $C_{\t,p; K}$ be the subspace of $K_{\t,p}$ consisting of all
constant functions. Similarly, set
 $$J_{\t,p}=(L_2(\el_2; t^{-\t}))_p^c\,\cap_p\,
 (L_2(\el_2; t^{1-\t}))_p^r\;,$$
and let $C_{\t,p; J}$ be the quotient space of $J_{\t,p}$ by the
subspace of mean zero functions. We refer to \cite{xu-int} for
more discussions on such spaces.

Note that $K_{\t, p}$ (resp. $J_{\t,p}$) is a quotient (resp. a
subspace) of $C_p\op_p R_p$. Thus both $C_{\t,p; K}$ and $C_{\t,p;
J}$ are quotients of subspace of $C_p\op_p R_p$.

\begin{thm}\label{Cq as a qs of Cp+Rp via real interpolation}
 Let $1\le p\le\8$, $\min(p, p')<q<\max(p,
p')$ and $\t$ be determined by $\frac{1}q=\frac{1-\t}p
+\frac{\t}{p'}$. Then $C_q=C_{\t,p; K}$ and $C_q=C_{\t,p; J}$
completely isomorphically with relevant constants majorized by
$2^{|1-\frac{2}p|}$.
\end{thm}

\pf We first show the complete inclusion $C_{\t,p; K}\subset C_q$
in the case of $p\ge 2$. Fix an $n\in\nat$ and let $x=\sum_kx_k\ot
e_k\in S_p^n[C_{\t,p; K}]$ of norm $<1$. (As before, we use
$(e_k)$ to denote the canonical basis of $\ell_2$.) Then there are
$f\in S_p^n[(L_2(\el_2; t^{-\t}))_p^c]$ and $g\in
S_p^n[(L_2(\el_2; t^{1-\t}))_p^r]$ such that $x=f(t)+g(t)$ for
almost all $t\in\real$ and
 $$\Big(\big\|f\big\|^{p}_{S_p^n[(L_2(\el_2;
t^{-\t}))_p^c]}
 \,+\,
 \big\|g\big\|^{p}_{S_p^n[(L_2(\el_2; t^{1-\t}))_p^r]}
 \Big)^{\frac{1}p}<1.$$
Using (\ref{description of Cp(Lp)}), we have
 \be
 \big\|f\big\|_{S_p^n[(L_2(\el_2;t^{-\t}))_p^c]}
 &=&\Big\|\int_0^\8\sum_k f_k(t)^*f_k(t)t^{-2\t}\,
 \frac{dt}{t}\Big\|^{\frac{1}2}_{S_{p/2}^n}\,,\\
 \big\|g\big\|_{S_p^n[(L_2(\el_2; t^{1-\t}))_p^r]}
 &=&\Big\|\int_0^\8\sum_k g_k(t)g_k(t)^*t^{2(1-\t)}\,\frac{dt}{t}
 \Big\|^{\frac{1}2}_{S_{p/2}^n}\,.
 \ee

To estimate the norm $\|x\|_{S_p^n[C_q]}$, we use again Lemma
\ref{description of norm in SpCq}. Let $\a$ (resp. $\b$) be a
positive unit element in $S^n_{2r\t^{-1}}$ (resp.
$S^n_{2r(1-\t)^{-1}}$). By perturbation, we can assume $\a$ and
$\b$ invertible. Let $A_0$ be the right multiplication operator on
$H=\el_2(S_2^n)$ by $\b^{\frac{1}{1-\t}}$: $(x_k)_k\mapsto
(x_k\b^{\frac{1}{1-\t}})_k$. Similarly, let $A_1$ be the left
multiplication operator by $\a^{\frac{1}\t}$. Then $A_0$ and $A_1$
are commuting invertible positive bounded operators on $H$. $A_i$
induces an equivalent norm on $H$: $\|x\|_{i}=\|A_i(x)\|$
($i=0,1$). Let $H_{i}$ be $H$ equipped with $\|\;\|_{i}$. Then
$(H_{0},\; H_{1})$ becomes a compatible couple of Hilbert spaces,
which can be identified as a compatible couple of weighted
$L_2$-spaces. Using the real interpolation for weighted
$L_2$-spaces (see \cite{bl}), we easily get the following:
 $$\|A_0^{1-\t}A_1^\t(x)\|=\|x\|_{(H_{0},\;
 H_{1})_{\t, 2; K}}\;, \quad x\in H.$$
Namely,
 $$\big(\sum_k \|\a x_k\b\|^2_2\big)^{\frac{1}{2}}
 =\|x\|_{(H_0,\; H_{1})_{\t, 2; K}}\;.$$
 The last interpolation norm is estimated as follows. First,
 \be
 \|f\|^2_{L_2(H_{0}; t^{-\t})}
 &=&\int_0^\8\sum_k \|f_k(t)\|_{H_0}^2\, t^{-2\t}\,\frac{dt}{t}\\
 &=&\int_0^\8\sum_k \Tr\big(\b^{\frac{1}{1-\t}}
 f_k(t)^*f_k(t)\b^{\frac{1}{1-\t}}\big)\, t^{-2\t}\,
  \frac{dt}{t}\\
 &=&\Tr\Big(\b^{\frac{2}{1-\t}}\int_0^\8\sum_k f_k(t)^*f_k(t)\, t^{-2\t}\,
  \frac{dt}{t}\Big)\\
 &\le& \Big\|\int_0^\8\sum_k
 f_k(t)^*f_k(t)t^{-2\t}\,\frac{dt}{t}\Big\|_{S^n_{p/2}}.
 \ee
Similarly,
 $$\|g\|^2_{L_2(H_{1}; t^{1-\t})}\le \Big\|\int_0^\8\sum_k
 g_k(t)g_k(t)^*t^{2(1-\t)}\,\frac{dt}{t}\Big\|_{S^n_{p/2}}.$$
Thus
 \be
 \big(\|f\|^2_{L_2(H_{0}; t^{-\t})}
 +\|g\|^2_{L_2(H_{1}; t^{1-\t})}\big)^{\frac{1}2}
 &\le& 2^{\frac{1}2 - \frac{1}p}\,
 \big(\|f\|^{p}_{S_p^n[(L_2(\el_2;
t^{-\t}))_p^c]}
 \,+\,
 \|g\|^{p}_{S_p^n [(L_2(\el_2; t^{1-\t}))_p^r]}
 \big)^{\frac{1}p}\\
 &<&2^{\frac{1}2 - \frac{1}p}
 \ee
Since $x=f(t)+g(t)$ for almost all $t\in\real_+$, we deduce that
 $$\|x\|_{(H_{0},\; H_{1})_{\t, 2; K}}
 <2^{\frac{1}2 - \frac{1}p}\;;$$
whence
 $$\big(\sum_k \|\a x_k\b\|^2_2\big)^{\frac{1}{2}}
 <2^{\frac{1}2 - \frac{1}p}\,,$$
which, together with Lemma \ref{description of norm in SpCq},
implies that $x\in S^n_p[C_q]$ and $\|x\|_{S_p[C_q]}\le
2^{\frac{1}2 - \frac{1}p}$. Therefore, we deduce the complete
inclusion $C_{\t,p; K}\subset C_q$ and that its c.b. norm is
majorized by $2^{\frac{1}2 - \frac{1}p}$ in the case of $p\ge2$.

To treat the case of $p\le 2$ we recall the fact that for a
Hilbert space $H$ we have
 $$(\overline{H}_{p'}^c)^*\cong H_{p}^c\cong H_{p'}^r\quad
 \mbox{completely isometrically}\,,$$
where $\overline H$ denotes the conjugate of $H$. Then using the
usual duality between sum and intersection spaces we obtain the
following complete isometric identifications:
 $$\big[(\,\overline{L_2(\el_2; t^{\t})}\,)_{p'}^c\,\cap_{p'}\,
 (\,\overline{L_2(\el_2; t^{-(1-\t)})}\,)_{p'}^r\big]^*
 \cong K_{\t, p}$$
and
 \be
 \big[(\,\overline{L_2(\el_2; t^{\t})}\,)_{p'}^c\,\cap_{p'}\,
 (\,\overline{L_2(\el_2; t^{-(1-\t)})}\,)_{p'}^r\big]^*
 &\cong & (L_2(\el_2; t^{-\t}))_{p'}^r\,+_p\,
 (L_2(\el_2; t^{(1-\t)}))_{p'}^c\\
 &\cong& (L_2(\el_2; t^{-(1-\t)}))_{p'}^c\,+_p\,
 (L_2(\el_2; t^{\t}))_{p'}^r\,.
 \ee
However, since $p\le2$,
 $$(L_2(\el_2; t^{-(1-\t)}))_{p'}^c\,+_p\,
 (L_2(\el_2; t^{\t}))_{p'}^r\subset K_{1-\t,p'}\quad
 \mbox{completely contractively}.$$
Therefore, we deduce that
 $$C_{\t, p; K}\subset C_{1-\t,p'; K}\quad
 \mbox{completely contractively}.$$
Applying the first part to $C_{1-\t,p'; K}$ (with $p'\ge2$), we
get that $C_{1-\t,p'; K}\subset C_q$ and the c.b. norm of this
inclusion is $\le 2^{\frac{1}2 -\frac{1}{p'}}$, and so  $C_{\t,p;
K}\subset C_q$ and is of c.b. norm $\le 2^{\frac{1}p
-\frac{1}{2}}$.

We have thus proved the complete inclusion $C_{\t,p; K}\subset
C_q$ for any $1\le p\le\8$ and its c.b. norm is majorized by
$2^{|\frac{1}2 -\frac{1}{p}|}$.

\medskip

To show the inverse inclusion we use the J-method. In a similar
way, we prove that $C_{\t,p; J}\subset C_q$ and is of c.b. norm
$\le 2^{|\frac{1}2 - \frac{1}p|}$. However, by the usual duality
between the K-  and J-methods, we easily check that
 $$\big[C_{\t,p; J}\big]^*=C_{1-\t,p' K}\quad\mbox{and}\quad
 \big[C_{\t,p; K}\big]^*=C_{1-\t,p' J}\quad \mbox{completely
 isometrically}.$$
It then follows that $C_q=C_{\t,p; K}$ and $C_q=C_{\t,p; J}$
completely isomorphically with relevant constants bounded by
$2^{|1 - \frac{2}p|}$. \cqd

\section{Noncommutative Khintchine inequalities}

This section contains the main ingredient for the embedding of the
quotients of subspace of $C_p\op_pR_p$ into noncommutative
$L_p$-spaces. It is the noncommutative Khintchine type
inequalities for Shlyakhtenko's generalized circular systems.

Given a complex Hilbert space $H$ we denote as usual by $\F(H)$
the associated free Fock space:
 $$\F(H)=\bigoplus_{n=0}^\8 H^{\otimes n},$$
where $H^{\otimes 0}={\com}\O$ with $\O$ a unit vector, called the
vacuum. Let $\ell(e)$ (resp. $\ell^*(e)$) denote the left creation
(resp. annihilation) operator associated with a vector $e\in H$.
Recall that $\ell^*(e)=(\ell(e))^*$.

Assume  $H$ is infinite dimensional and separable. Fix an
orthonormal basis $\{e_{\pm k}\}_{k\ge 1}$ in $H$. We will also
fix a sequence $\{\l_k\}_{k\ge 1}$ of positive numbers. Let
 $$g_k=\el(e_{k}) + \sqrt{\l_k}\,\el^*(e_{-k}).$$
$(g_k)_{k\ge 1}$ is a generalized circular system in
Shlyakhtenko's sense \cite{shlya-quasifree}. Let $\Ga$ be the von
Neumann algebra on $\F(\H)$ generated by the $g_k$. Let $\rho$ be
the vector state on $\Ga$ determined by $\O$. By
\cite{shlya-quasifree}, $\rho$ is faithful on $\Ga$. By the
identification of $L_1(\Ga)$ with $\Ga_*$, $\rho$ is a positive
unit element in $L_1(\Ga)$, so for any $1\le p\le\8$, $\rho^{1/p}$
is a positive unit element in $L_p(\Ga)$, and thus
$g_k\,\rho^{1/p}\in L_p(\Ga)$ for any $k$.

\medskip

The following is the noncommutative Khintchine type inequalities
for generalized circular systems from \cite{xu-gro}. See also
\cite{jx-free} for this kind of inequalities in a more general
setting. Note that the case $p=\8$ is already contained in
\cite{pisshlyak}.

\begin{thm}\label{khintchine}
 Keep the notations above. Let $M$ be a von
Neumann algebra  and $(x_n)$  a finite sequence in $L_p(M)$.
 \begin{enumerate}[{\rm i)}]
 \item  If $p\ge 2$,
 \be
 &&\max\Big\{
 \big\|\big(\sum_k x_k^*x_k\big)^{\frac{1}{2}}\big\|_{L_p(M)}\;,\
 \big\|\big(\sum_k\l_k^{1-\frac{2}{p}}\, x_kx_k^*
 \big)^{\frac{1}{2}}\big\|_{L_p(M)} \Big\}\\
 &&\hskip 2.5cm \le\big\|\sum_kx_k\ot g_k\rho^{\frac{1}p}
 \big\|_{L_p(M\bar\ot\Ga)}\le \\
 &&B_p\max\Big\{\big\|\big(\sum_k x_k^*x_k
 \big)^{\frac{1}{2}}\big\|_{L_p(M)}\;,\
 \big\|\big(\sum_k\l_k^{1-\frac{2}{p}}\, x_kx_k^*
 \big)^{\frac{1}{2}}\big\|_{L_p(M)} \Big\}.
 \ee
 \item  If $p<2$,
 \be
 && A_p^{-1} \inf\Big\{\big\|\big(\sum_k a_k^*a_k
 \big)^{\frac{1}{2}}\big\|_{L_p(M)} +
 \big\|\big(\sum_k\l_k^{1-\frac{2}{p}}\, b_kb_k^*
 \big)^{\frac{1}{2}}\big\|_{L_p(M)}\Big\}\\
 &&\hskip 2.5cm \le\big\|\sum_kx_k\ot g_k\rho^{\frac{1}p}
 \big\|_{L_p(M\bar\ot\Ga)}\le\\
 &&\inf\Big\{\big\|\big(\sum_ka_k^*a_k\big)^{\frac{1}{2}}
 \big\|_{L_p(M)} +
 \big\|\big(\sum_k\l_k^{1-\frac{2}{p}}\, b_kb_k^*
 \big)^{\frac{1}{2}}\big\|_{L_p(M)}\Big\},
 \ee
where the infimum runs over all decompositions $x_k=a_k+b_k$ in
$L_p(M)$. The two positive constants $A_p$ and $B_p$ depend only
on $p$ and can be controlled by a universal constant.
  \item  Let $G_{p}$ be the closed subspace of
$L_p(\Ga)$ generated  by $\{g_k\rho^{1/p}\}_{k\ge1}$. Then $G_{p}$
is completely complemented in  $L_p(\Ga)$ with constant
 $\le 2^{|1-\frac{2}{p}|}$.
 \end{enumerate}
\end{thm}

By \cite{shlya-quasifree} the algebra $\Ga$ generated by the $g_k$
is a type III$_\l$ factor ($0<\l\le1$) (except the case where all
$\l_k$ are equal to 1, which is the usual case of Voiculsecu's
circular systems). $\Ga$ is never hyperfinite. Recall that $\Ga$
is the free analogue of the classical Araki-Woods quasi-free CAR
factors. The latter factors are \underbar {hyperfinite} type
III$_\l$. Thus it is natural to ask whether similar Khintchine
inequalities also hold in the classical quasi-free case. This is
indeed true for $1<p<\8$.  To state this result, we need to
introduce the Fermionic analogue of the system $(g_k)$.

\medskip

Again consider a complex Hilbert space $H$ as before. We denote by
$\F_{-1}(H)$ the associated antisymmetric Fock space:
 $$\F_{-1}(H)=\bigoplus_{n=0}^\8 \Lambda^n(H),$$
where $\Lambda^n(H)$ stands for the $n$-fold antisymmetric tensor
product of $H$  (with $\Lambda^0(H)={\com}\O$). Let $c(e)$ (resp.
$c^*(e)$) denote the creation (resp. annihilation) operator
associated with a vector $e\in H$. Recall that $c^*(e)=(c(e))^*$.

Let $\{e_{\pm k}\}_{k\ge 1}$ and $\{\l_k\}_{k\ge 1}$  be as
before. Put
 $$f_k=c(e_{k}) + \sqrt{\l_k}\,c^*(e_{-k}).$$
Then $(f_k)_{k\ge 1}$ is a CAR sequence. Let $\Ga_{-1}$ be the von
Neumann algebra on $\F_{-1}(H)$ generated by the $f_k$ and
$\rho_{-1}$ the vector state on $\Ga_{-1}$ determined by $\O$.
$\rho_{-1}$ is a quasi-free state, which is faithful on
$\Ga_{-1}$.

The following is the Fermionic analogue of
Theorem~\ref{khintchine}.

\begin{thm}\label{khintchine-fermi}
 Let $M$ be a von Neumann algebra, $1<p<\8$
 and $(x_k)\subset L_p(M)$ a finite sequence.
 \begin{enumerate}[\rm i)]
 \item  If $2\le p<\8$,
 \be
 &&\max\Big\{\big\|\big(\sum_k\, x_k^*x_k\big)^{\frac{1}2}
 \big\|_{L_p(M)}\ ,\
 \big\|\big(\sum_k\l_k^{1-{\frac{2}p}}\,x_kx^*_k\big)^{\frac{1}2}
 \big\|_{L_p(M)}\Big\}\\
 &&\hskip 2.2cm\le\big\|\sum_kx_k\ot
 f_k\rho_{-1}^{\frac{1}p}\big\|_{L_p(M\bar\ot\Ga_{-1})}\le\\
 &&B_p'
 \max\Big\{\big\|\big(\sum_k\, x_k^*x_k\big)^{\frac{1}2}
 \big\|_{L_p(M)}\ ,\
 \big\|\big(\sum_k\l_k^{1-{\frac{2}p}}\, x_kx^*_k\big)^{\frac{1}2}
 \big\|_{L_p(M)}\Big\}\ .
 \ee
 \item  If $1<p<2$,
 \be
 &&A_p' \inf\Big\{
 \big\|\big(\sum_k\, a_k^*a_k\big)^{\frac{1}2}\big\|_{L_p(M)}
 +\big\|\big(\sum_k\l_k^{1-\frac{2}p}\, b_kb^*_k\big)^{\frac{1}2}
 \big\|_{L_p(M)}\Big\}\\
 &&\hskip 2.2cm\le \big\|\sum_k x_k\ot
 f_k\rho_{-1}^{\frac{1}p}\big\|_{L_p(M\bar\ot\Ga_{-1})}\le\\
 &&\inf\Big\{\big\|\big(\sum_k\, a_k^*a_k\big)^{\frac{1}2}
 \big\|_{L_p(M)}
 +\big\|\big(\sum_k\l_k^{1-\frac{2}p}\, b_kb^*_k\big)^{\frac{1}2}
 \big\|_{L_p(M)}\Big\}\,,
 \ee
where the infimum runs over all decompositions $x_k=a_k+b_k$ with
$a_k, b_k\in L_p(M)$. Here the positive constants $A_p'$ and $
B_p'$ depend only on $p$.
\item  Let $F_p$ be the closed subspace of $L_p(\Ga_{-1})$
generated  by $\{f_k\rho_{-1}^{\frac{1}p}\}_k$. Then $F_p$ is
completely complemented in $L_p(\Ga_{-1})$.
\end{enumerate}
\end{thm}

\pf Part i) follows from the noncommutative Burkholder inequality
from \cite{jx-burk} (see also \cite{jx-rosenthal}). Parts ii) and
iii) are proved in a similar way as their free counterparts in
Theorem \ref{khintchine} (see \cite{xu-gro}). We omit the details.
\cqd

\medskip\n{\bf Remark.}
Very recently, Junge \cite{ju-khintchine fermi} proved that
Theorem \ref{khintchine-fermi}, iii) holds for $p=1$ too. However,
it seems still unknown whether the constant $A'_p$ there can be
controlled by an absolute constant.

\section{Embeddings}

This section is devoted to the problem of embedding $C_q$ (and
more generally the quotients of subspace of $C_p\op_pR_p$) into
noncommutative $L_p$-spaces for $1\le p<q\le 2$.

\begin{thm}\label{embed qs}
 Let $H$ and $K$ be two Hilbert spaces.
Let $1\le p<2$ and $X$ be a quotient of subspace of $H_p^c\op_p
K_p^r$. Then there is a QWEP von Neumann algebra $M$ such that $X$
embeds completely isomorphically into $L_p(M)$. Moreover, the
relevant constant can be majorized by a universal constant.\\
 \indent If $p>1$, $M$ can be chosen hyperfinite.
\end{thm}

 \pf For simplicity, we assume $H=K=\el_2$. Let $X$ be a quotient of
subspace of $C_p\op_p R_p$. Then $Y=X^*$ is a quotient of subspace
of $C_{p'}\op_{p'} R_{p'}$, so there is a subspace $E\subset
C_{p'}\op_{p'} R_{p'}$ such that $Y$ is a quotient of $E$. By
Corollary \ref{graphter}, $E$ is decomposed into a direct sum of
three spaces: a $p'$-column space, a $p'$-row space and the graph
of an injective closed densely defined operator with dense range
from a $p'$-column space into a $p'$-row space. Clearly, we need
only to consider the third term. Thus without loss of generality,
assume $E=G(\D)\subset C_{p'}\op_{p'} R_{p'}$, where $\D:
\ell_2\to \ell_2$ is a positive injective closed densely defined
operator with dense range (see Remark \ref{graphbis}). By
discretization, we may as well assume $\D$ diagonal:
$\D(e_k)=\d_k\,e_k$ for all $k$, where $(e_k)$ is the canonical
basis of $\ell_2$ and $(\d_k)$ a sequence of positive numbers.
Note that $E$ is simply the diagonal subspace of $ C_{p'}\op_{p'}
R_{p'}$ but with the second space weighted by the $\d_k$. Passing
to duals, we deduce that $E^*$ is the quotient of $C_p\op_p R_p$
by the subspace $\{(x, -x): x\in\ell_2\}$, i.e. $E^*=C_p+_p R_p$.

Now we are in a position to apply Theorem~\ref{khintchine}. Let
$\l_k$ be such that $\d_k=\l_k^{1/2-1/p}$. Then by the above
discussion and Theorem~\ref{khintchine}, ii), we see that $E^*$ is
completely isomorphic to $G_p$ in Theorem~\ref{khintchine}.
Therefore, $X=X^{**}=Y^*$ is completely isomorphic to a subspace
of $G_p\subset L_p(\Ga)$. By \cite{pisshlyak}, $\Ga$ is QWEP.

If $p>1$, instead of Theorem~\ref{khintchine}, we can use
Theorem~\ref{khintchine-fermi}. Then $E^*$ is completely
isomorphic to $F_p$ in Theorem~\ref{khintchine-fermi}.\cqd

\medskip

\begin{rk}\label{embed graph}
 i) By Theorem \ref{khintchine}, the space $G_p$ is
completely complemented in $L_p(\Ga)$. Thus the proof above shows
that the graph of an operator from $C_p$ into $R_p$ is completely
isomorphic to a completely complemented subspace of $L_p(\Ga)$.
This fact is useful for  applications.

ii) By \cite{ju-khintchine fermi}, the von Neumann algebra $M$ in
Theorem \ref{embed qs} can be chosen to be hyperfinite in the case
of $p=1$ too.
\end{rk}

Let $QS(C_p, R_p)$ denote the family of all operator spaces which
are completely isomorphic to quotients of subspace of $H_p^c\op_p
K_p^r$ for some Hilbert spaces $H$ and $K$. It is clear that if
$X$ belongs to $QS(C_p, R_p)$, then so does $X^*$. It was proved
in \cite{xu-gro} that if both $X$ and $X^*$ embed completely
isomorphically into some $L_p(M)$ ($1<p<2$), then $X\in QS(C_p,
R_p)$. This together with the previous theorem implies the
following

\begin{cor}
 Let $1<p<2$ and $X$ be an operator space. Then $X\in QS(C_p,R_p)$
iff both $X$ and $X^*$ embed completely isomorphically into a
noncommutative $L_p$-space.
\end{cor}

\n{\bf Remark.} It was shown in \cite{pisshlyak} that if $X$ is
exact and both $X$ and $X^*$ completely embed into an $L_1$, then
$X\in QS(C_1, R_1)$.

\begin{thm}\label{embed Cq}
 Let $1\le p<q\le 2$. Then there is a QWEP
type III factor $M$ such that $C_q$ embeds completely
isomorphically into $L_p(M)$.\\
\indent The same assertion holds for $R_q$ too.
\end{thm}

\pf By Theorem~\ref{Cq as a qs of Cp+Rp}, $C_q$ is a quotient of
subspace of $C_p\op_pR_p$. Thus Theorem~\ref{embed qs} implies
that $C_q$ embeds into an $L_p(M)$. Let us show that $M$ can be
taken to be a type III factor.  Let $E$ be the space introduced
during the proof of Theorem~\ref{Cq as a qs of Cp+Rp}. Since an
analytic function on the strip $S$ is uniquely determined by its
values on $\partial_0$ (or $\partial_1$), we see that $E$ is a
graph. Thus by Remark \ref{embed graph}, $E$ embeds into
$L_p(\Ga)$. By \cite{shlya-quasifree}, $\Ga$ is a type III
factor.\cqd

\begin{rk}
 The factor $\Ga$ above is of type III$_\l$ for some $\l\in(0,1]$.
Using Theorem \ref{Cq as a qs of Cp+Rp via real interpolation} and
the precise form of the space $K_{\t, p}$ introduced there, one
easily shows that for any $\l\in(0,1]$, $\Ga$ can be chosen to be
of type III$_\l$. This follows from an appropriate  discretization
of $K_{\t, p}$ and the classification theorem in
\cite{shlya-quasifree}.
\end{rk}

\medskip

Pisier \cite{pis-III} proved that $C_q$ ($1<q\le 2$) cannot embed
completely isomorphically  into the predual of a semifinite von
Neumann algebra. We now extend this result to our general setting.

\begin{thm}\label{embed semifinite}
 Let $1\le p<q\le 2$. Then neither $C_q$ nor
$R_q$ embeds completely isomorphically into a noncommutative
$L_p(M)$ with a semifinite von Neumann algebra $M$.
\end{thm}

Junge \cite{ju-emb} proved that with the same $p, q$ as above, the
Schatten class $S_q$ embeds \underbar{isomorphically} into
$L_p(\R)$ at the Banach space level, where $\R$ is the hyperfinite
factor of type II$_1$. It was left open there whether one can do
this at the operator space level. Theorem~\ref{embed semifinite}
clearly implies that this is impossible. Thus we have the

\begin{cor}
 Let $1\le p<q\le 2$. Then  $S_q$ cannot embed completely
isomorphically into a semifinite noncommutative $L_p$-space.
 \end{cor}

\n{\bf Remark.} It is still unknown whether $S_q$ embeds
completely isomorphically into a type III noncommutative $L_p$.

\medskip

For the proof of Theorem~\ref{embed semifinite} we will need the
following factorization of c.b. maps from $L_p(M)$ to $C_q$.

\begin{lem}\label{factorization}
 Let $M$ be a von Neumann algebra  and $H$ a Hilbert space.
Let $2\le p\le\8,\; 0\le\t\le1$ and
$\frac{1}{q}=\frac{1-\t}{p}+\frac{\t}{p'}$. Then for any map $u:
L_p(M)\to H_q^c$ the following assertions are equivalent:
 \begin{enumerate}[\rm i)]
\item  $u$ is c.b..
\item  There are two positive unit functionals $f, g\in
\big(L_{p/2}(M)\big)^*$ such that
 \beq\label{factorization1}
 \|u(x)\|\le c\, \big(f(x^*x)\big)^{(1-\t)/2}
 \big(g(xx^*)\big)^{\t/2}\ ,\quad x\in L_p(M).
 \eeq
\end{enumerate}
Moreover, if $c$ denotes the best constant in
(\ref{factorization1}), then
 $$c_{p,\,\t}^{-1}\;c\le\|u\|_{cb}\le  c,$$
where $c_{p,\,\t}$ is a positive constant depending only on $p$
and $\t$, which can be controlled by an absolute constant.
\end{lem}

This is proved in \cite{pis-Rp} for $p=\8$ and \cite{xu-gro} for
$p<\8$. Moreover, the same result holds for any $u$ defined on a
subspace $E$ of $L_p(M)$ ($E$ must be then assumed exact in the
case $p=\8$).

\medskip

We will also need the following simple description of c.b. maps
between $C_p$ and $C_q$. This is probably known to experts.

\begin{lem}\label{cb CpCq}
 Let $1\le p, q\le\8$ with $p\neq q$. Then
 $$CB(C_p,\; C_q)=S_{\frac{2pq}{|p-q|}}\quad\mbox{isometrically}.$$
\end{lem}

\pf We first consider the case where $p=\8$. We have
 $$CB(C,\; C_q)=R\ot_{\min}C_q
 =\big(R\ot_{\min}C,\ R\ot_{\min}R\big)_{1/q}
 =\big(S_\8,\ S_2\big)_{1/q}=S_{2q}.$$

Now assume $p>q$.  Let $u\in CB(C_p, C_q)$. For any $v\in CB(C,
C_p)$, $uv\in CB(C, C_q)$ and
 $$\|uv\|_{CB(C,\; C_q)}\le \|u\|_{CB(C_p,\; C_q)}
 \|v\|_{CB(C,\; C_p)}.$$
Thus
 $$\|uv\|_{S_{2q}}\le \|u\|_{CB(C_p,\; C_q)}
 \|v\|_{S_{2p}}.$$
Taking the surepmum over all $v$ in the unit ball of $S_{2p}$, we
get
 $$\|u\|_{S_{\frac{2pq}{p-q}}}\le \|u\|_{CB(C_p,\; C_q)}.$$
To prove the converse inequality we use interpolation. Let
$\t=q/p$. Then
 $$C_p=(C,\ C_q)_\t.$$
Thus by the case above and  interpolation
 \be
 S_{\frac{2pq}{p-q}}
 &=&(S_{2q},\ B(\ell_2))_\t
 =\big(CB(C,\; C_q),\ CB(C_q,\; C_q)\big)_\t\\
 &\subset& CB\big((C,\ C_q)_\t,\ C_q\big)
 =CB(C_p,\; C_q),
 \ee
where the inclusion is completely contractive. Therefore, we have
obtained
 $$CB(C_p,\; C_q)=S_{\frac{2pq}{p-q}}\quad\mbox{for}\ p> q.$$
Passing to adjoints, we deduce the identity in the case of
$p<q$.\cqd

\medskip\n{\bf Remark.} The lemma above immediately yields the
following isometric identities:
 $$CB(R_p,\; R_q)=S_{\frac{2pq}{|p-q|}}\quad\mbox{and}\quad
  CB(R_p,\; C_q)=S_{\frac{2pq}{|p+q-pq|}}=CB(C_p,\; R_q).$$

\medskip

\n{\it Proof of Theorem~\ref{embed semifinite}.} Assume $C_q$ is
completely isomorphic to a subspace of $L_p(M)$, where $M$ is a
semifinite von Neumann algebra equipped with a normal semifinite
faithful trace $\tau$. For simplicity, assume $C_q$ itself is a
subspace of $L_p(M)$. Let $J: C_q\to L_p(M)$ be the inclusion map.
Then $\displaystyle u\mathop{=}^{\rm def}J^*: L_{p'}(M)\to C_{q'}$
is completely contractive. Recall that $1\le p<q\le 2$. Thus $p'$
and $q'$ satisfy the assumption of Lemma~\ref{factorization}.
Therefore, by that lemma, there are positive functionals $f$ and
$g$ in the unit ball of $\big(L_{p'/2}(M)\big)^*$ such that
(\ref{factorization1}) holds. Replacing $f$ and $g$ by $f+g$, we
can assume $f=g$. Thus we have
 $$\|u(x)\|\le c\, \big(f(x^*x)\big)^{(1-\t)/2}
 \big(f(xx^*)\big)^{\t/2}\ ,\quad x\in L_{p'}(M).$$
If $p>1$, $f$ can be viewed as a positive unit operator in
$L_r(M)$, where $r$ is the conjugate index of $p'/2$. On the other
hand, if $p=1$, $f$ can be taken as a normal state on $M$ for $u$
is normal (see \cite{pis-III} for more details). Thus in this
case, $f$ is again a positive unit operator in $L_1(M)$.
Therefore, the above inequality can be rewritten as
 \beq\label{factorization2}
 \|u(x)\|\le c\, \big(\tau(f x^*x)\big)^{(1-\t)/2}
 \big(\tau(f xx^*)\big)^{\t/2}\ ,\quad x\in L_{p'}(M).
 \eeq
As in \cite{pis-III}, we are going to reduce $M$ to a finite von
Neumann algebra. To this end, given $\l>1$ let $e_\l$ be the
spectral projection of $f$ corresponding to the interval
$(\l^{-1},\l)$. Since $f\in L_r(M)$ with $r<\8$, $\tau(e_\l)<\8$
and $e_\l\to 1$ strongly as $\l\to\8$. Decompose $J$ as follows:
 $$J=R_{e_\l}\,L_{e_\l}J+R_{e_\l}\,L_{e_\l^\perp}J+
 R_{e_\l^\perp}\,J\buildrel{\rm def}\over
 =J_1 +J_2+J_3,$$
where $L_a$ (resp. $R_a$) denotes the left (resp. right)
multiplication by $a$. Consequently,
 $$u=uR_{e_\l}\,L_{e_\l}+uL_{e_\l}\,R_{e_\l^\perp}+
 uL_{e_\l^\perp}\mathop=^{\rm def}\;u_1 +u_2+u_3.$$
We show that $\|u_2\|_{cb}$ and $\|u_3\|_{cb}$ become very small
for $\l$ sufficiently large. By (\ref{factorization2}), for any
$x\in L_{p'}(M)$
 \be
 \|u_2(x)\|^2
 &=&\|u(e_\l xe_\l^\perp)\|^2
 \le c^2\, \big(\tau(f e_\l^\perp x^*e_\l xe_\l^\perp)\big)^{1-\t}
 \big(\tau(f e_\l xe_\l^\perp x^*e_\l)\big)^{\t}\\
 &\le& c^2\, \big(\tau(fe_\l^\perp x^*x)\big)^{1-\t}
 \big(\tau(f  x x^* )\big)^{\t}\\
 &=& c^2\,\|fe_\l^\perp\|_r^{1-\t} \big(\tau(g x^*x)\big)^{1-\t}
 \big(\tau(f  x x^* )\big)^{\t}\,,
 \ee
where $g=\|fe_\l^\perp\|_r^{-1}fe_\l^\perp $ is a positive unit
operator in $L_r(M)$. Thus applying Lemma~\ref{factorization} to
$u_2$, we get
 $$\|u_2\|_{cb}\le c'\|fe_\l^\perp\|_r^{(1-\t)/2}.$$
Since $e_\l^\perp\to 0$ strongly as $\l\to\8$,
$\|fe_\l^\perp\|_r\to 0$. Therefore, $\|u_2\|_{cb}$ is small for
large $\l$. Similarly, we get the the same assertion for $u_3$.
Passing to adjoints, we deduce that  $\|J_2\|_{cb}$ and
$\|J_3\|_{cb}$ are small for large $\l$. Consequently, $J_1$ is a
perturbation of $J$, and so $J_1$ is a complete isomorphism (if
$\l$ is large enough). However, $J_1$ takes values in $L_{p'}(e_\l
Me_\l)$ and the restriction of $\tau$ to $e_\l Me_\l$ is finite.
On the other hand, $fe_\l$ is bounded by $\l$. Therefore,
replacing $J$ by $J_1$ if necessary, we are reduced to the case
where $\tau$ is finite, say, normalized, and where $f$ is a
positive bounded operator in $M$, say, the identity of $M$. Then
(\ref{factorization2}) becomes,
 \beq\label{factorization3}
 \|u(x)\| \le c\, \big(\tau(x^*x)\big)^{1/2},
 \quad x\in L_{p'}(M).
 \eeq

The rest of the proof is to show that $J$ factors through a
$p$-column space. For this  we consider the scalar product $\la y,
x\ra=\tau(y^*x)$ on $L_{p'}(M)$. Let $H$ be the completion of
$L_{p'}(M)$ with respect to this scalar product. It is clear that
the identity on $L_{p'}(M)$ induces a contractive inclusion of
$L_{p'}(M)$ into $H$, denoted by $\iota$. On the other hand,
(\ref{factorization3}) implies that there is a bounded operator
$T: H\to C_{q'}$ such that $\|T\|\le c$ and $u=T\iota$. We now
equip $H$  with the operator space structure of $H_{p'}^c$. Then
$\iota$ is completely contractive. This was already proved in
\cite{xu-gro}. We include a proof for completeness. By
\cite[Lemma~1.7]{pis-ast}, it suffices to prove that
$I_{S_{p'}}\ot\iota$ extends to a contraction from
$S_{p'}[L_{p'}(M)]$ to $S_{p'}[H_{p'}^c]$. Let $x=(x_{ij})\in
S_{p'}[L_{p'}(M)]$ be a finite matrix. Then
 $$\|I_{S_{p'}}\ot\iota (x)\|^2_{S_{p'}[H_{p'}^c]}
 =\big\|\big(\tau(\sum_kx_{ki}^*x_{kj})\big)_{ij}
 \big\|_{S_{p'/2}}.$$
Let $\a$ be a positive finite matrix in the unit ball of $S_{t}$
($t$ being the conjugate index of $p'/2$). Then (with $\Tr$ the
usual trace on $B(\ell_2)$)
 \be
 \sum_{i,j}\a_{ij}\tau(\sum_kx_{ki}^*x_{kj})
 &=&\Tr\ot\tau\big((\a\ot 1)x^*x)\\
 &\le& \|\a \ot 1\|_{S_{t}[L_{t}(M)]}\,
 \|x^*x\|_{S_{p'/2}[L_{p'/2}(M)]}
 \le \|x\|^2_{S_{p'}[L_{p'}(M)]}.
 \ee
Then taking the supremum over $\a$, we obtain
 $$\|I_{S_{p'}}\ot\iota (x)\|^2_{S_{p'}[H_{p'}^c]}
 \le \|x\|^2_{S_{p'}[L_{p'}(M)]}.$$
Therefore, $\iota $ is completely contractive.

Returning back to the factorization $u=T\iota$ and passing to
adjoints, we deduce that
 $\displaystyle J=\iota^*T^*\,\mathop=^{\rm def}\,vw$.
Hence the identity of $C_q$ factors through $\bar H_p^c$ with
$\|v\|_{cb}\le1$. Restricted to a subspace of $\bar H_p^c$ if
necessary, $v$ can be assumed to have values in $J(C_q)=C_q$. Then
$H$ must be separable and infinite dimensional. By Lemma \ref{cb
CpCq},  $v\in S_{\frac{2pq}{q-p}}$. Consequently, $J\in
S_{\frac{2pq}{q-p}}$, a contradiction, which achieves  the proof
of Theorem~\ref{embed semifinite}.\cqd

\bigskip

\n{\bf Acknowledgements.} We are very grateful to Eric Ricard for
fruitful discussions.
\bigskip


\end{document}